# THE BI-POISSON PROCESS: A QUADRATIC HARNESS

By Włodzimierz Bryc ,[1] Wojciech Matysiak
and Jacek Wesołowski

*University of Cincinnati, Warsaw University of Technology
and Warsaw University of Technology*

This paper is a continuation of our previous research on quadratic harnesses, that is, processes with linear regressions and quadratic conditional variances. Our main result is a construction of a Markov process from given orthogonal and martingale polynomials. The construction uses a two-parameter extension of the Al-Salam–Chihara polynomials and a relation between these polynomials for different values of parameters.

**1. Introduction.** The family of stochastic processes with uncorrelated increments, linear regressions and quadratic conditional variances includes the Wiener, Poisson, Gamma and Pascal processes; for details and additional references, see [20], Theorem 1. It includes their free counterparts which are Markov processes whose bivariate distributions match the bivariate distributions of free Brownian motion and free Poisson process; for details and additional references, see [4], [10], Theorem 4.3, and [5], Proposition 3.4. It also includes the classical version of the $q$-Brownian motion, see [6], and the classical version of the $q$-Poisson process, see [1]. All of these examples share the property that the conditional variance with respect to the past sigma field is constant.

In this paper we study the three-parameter family of stochastic processes with linear rather than constant conditional variances with respect to the past or the future sigma fields. The bi-Poisson process with parameters $(\eta, \theta, q)$ is the square-integrable stochastic process $(X_t)$ such that for all

Received May 2006; revised February 2007.
[1]Supported in part by NSF Grants INT-03-32062, DMS-05-04198 and by the C. P. Taft Memorial Fund.

$0 \leq s < t < u$ we have

$$(1) \qquad \mathbb{E}(X_t) = 0, \qquad \mathbb{E}(X_t X_s) = \min\{t, s\}.$$

$$(2) \qquad \mathbb{E}[X_t | \mathcal{F}_{s,u}] = \frac{u-t}{u-s} X_s + \frac{t-s}{u-s} X_u$$

$$
\begin{aligned}
(3) \quad \mathrm{Var}[X_t | \mathcal{F}_{s,u}] = \frac{(u-t)(t-s)}{u-qs} \Bigg( & 1 + \eta \frac{uX_s - sX_u}{u-s} \\
& + \theta \frac{X_u - X_s}{u-s} \\
& - (1-q) \frac{(uX_s - sX_u)(X_u - X_s)}{(u-s)^2} \Bigg),
\end{aligned}
$$

where

$$\mathcal{F}_{s,u} = \sigma\{X_t : 0 \leq t \leq s \text{ or } t \geq u\}.$$

Thus in the terminology of [8], the bi-Poisson process is a quadratic harness with parameters $\tau = \sigma = 0$.

We use the standard $q$-notation

$$[n]_q = 1 + q + \cdots + q^{n-1},$$

$$[n]_q! = [1]_q [2]_q \cdots [n]_q,$$

$$\begin{bmatrix} n \\ k \end{bmatrix}_q = \frac{[n]_q!}{[n-k]_q! [k]_q!},$$

with the usual conventions $[0]_q = 0, [0]_q! = 1$.

In [8], Proposition 4.13, we give the following uniqueness result.

PROPOSITION 1.1.   *Fix* $-1 \leq q \leq 1$. *If polynomials* $\{p_n(x;t)\}$ *given by*

$$
\begin{aligned}
(4) \quad x p_n(x;t) = {} & p_{n+1}(x;t) + (\theta + t\eta)[n]_q p_n(x;t) \\
& + t(1 + \eta\theta[n-1]_q)[n]_q p_{n-1}(x;t),
\end{aligned}
$$

$n \geq 0$, *with* $p_{-1} = 0, p_0 = 1$, *are orthogonal martingale polynomials for a Markov process* $(X_t)$, *then* $(X_t)$ *is a bi-Poisson process, that is,* (2) *and* (3) *hold, and* $1 + \eta\theta \geq \max\{q, 0\}$. *Moreover, such a Markov process* $(X_t)$ *is determined uniquely.*

The goal of this paper is to construct the three-parameter family of bi-Poisson processes within the full range of the parameters.

THEOREM 1.2.   *For* $-1 \leq q \leq 1$ *and* $1 + \eta\theta \geq \max\{q, 0\}$ *there exists a Markov process* $(X_t)$ *with* $X_0 = 0$ *and transition probabilities such that* $\{p_n(x;t)\}$ *are martingale and orthogonal polynomials for* $(X_t)$.

Combining Proposition 1.1 with Theorem 1.2 we get the following.

COROLLARY 1.3. *For every* $-1 \leq q \leq 1$, $1 + \eta\theta \geq \max\{q, 0\}$*, there exists a unique bi-Poisson process* $(X_t)$ *with parameters* $(\eta, \theta, q)$*. Moreover,* $(X_t)$ *is Markov.*

Polynomials (4) are the so-called Al-Salam–Chihara polynomials and have been studied by many authors; see [13], page 377, [14], Section 3.8. A natural plan of constructing the Markov process corresponding to a prescribed family of orthogonal polynomials relies on the Mehler–Poisson formula for the polynomials; such a plan was followed successfully in [7] where the polynomials were simple enough with $\eta = \theta = 0$. However, the Mehler–Poisson formula for the more general polynomials (4) is not available. We therefore follow the approach from [10], where we define the Markov transition probabilities of $(X_t)$ using auxiliary orthogonal polynomials $\{Q_n\}$, which were guessed by a trial and error method. This approach was then tested on the bi-Poisson process with parameter $q = 0$ in [12], where polynomials $\{p_n\}$ and also $\{Q_n\}$ satisfy a constant coefficients recurrence. The difficulty of the general bi-Poisson case lies in the fact that the corresponding orthogonal polynomials $\{Q_n\}$ satisfy a more complex three step recurrence which falls outside of the classes studied in the literature. The proof of the Markov property relies on nontrivial identities for connection coefficients between these orthogonal polynomials; the form of the identities deviates significantly from their predecessors [10], Lemma 3.1, [12], Proposition 2.2, and [9], Theorem 1; the latter identity expressed the Al-Salam–Chihara polynomials in terms of the continuous $q$-Hermite polynomials and provided the initial breakthrough that made our constructions feasible.

**Note added late.** After the first version of this paper was written, we learned about the so-called "projection formulas" for the Askey–Wilson and related polynomials in papers by [17, 18, 19]. From the results in these papers one can deduce that there is a family of probability measures $\{P_{s,t}(x, dy)\}$ such that the integral representation

$$p_n(x; t) = \int p_n(y; u) P_{t,u}(x, dy)$$

holds for all $t \leq u$, all $u$ from the (possibly empty) interval

$$(5) \qquad \left( \frac{\theta^2}{\eta\theta + 1 - q}, \frac{\eta\theta + 1 - q}{\eta^2} \right),$$

and for all $x$ from the interval (15). One can further check that the Chapman–Kolmogorov equations hold for $s \leq t$ from the interval (5), so measures $\{P_{s,t}(x, dy)\}$ are the transition probabilities of the Markov process from Theorem 1.2. Thus one can deduce the existence of the process $(X_t)$ for $t$ in (5).

The proofs in [17, 18, 19] rely on infinite product identities for hypergeometric functions and do not cover the full range $0 < t < \infty$. On the other hand, Theorem 1.2 does not generalize [17, 18, 19], as in our setting two of their parameters vanish.

**2. Orthogonal polynomials and the bi-Poisson process.** We will determine the transition probabilities of the general bi-Poisson process in implicit way as the orthogonality measures of the auxiliary six-parameter family of orthogonal polynomials $\{Q_n(y; x, t, s) : n \geq -1\}$ which satisfy the three step recurrence

(6)
$$
\begin{aligned}
yQ_n&(y; x, t, s) \\
&= Q_{n+1}(y; x, t, s) + \mathcal{A}_n(x, t, s)Q_n(y; x, t, s) \\
&\quad + \mathcal{B}_n(x, t, s)Q_{n-1}(y; x, t, s),
\end{aligned}
$$

with $Q_{-1} \equiv 0$, $Q_0 \equiv 1$ and with the coefficients

(7)
$$
\mathcal{A}_n(x, t, s) = q^n x + [n]_q(t\eta + \theta - [2]_q q^{n-1}s\eta),
$$

(8)
$$
\mathcal{B}_n(x, t, s) = [n]_q(t - sq^{n-1})\{1 + \eta xq^{n-1} + [n-1]_q\eta(\theta - s\eta q^{n-1})\},
$$

where $n \geq 0$. [For $n = 0$, the formulas should be interpreted as $\mathcal{A}_0(x, t, s) = x$, $\mathcal{B}_0(x, t, s) = 0$.]

In order for such a definition of transition probabilities to make sense, we will need to establish a number of properties of these polynomials. First, we need to know that the orthogonality measure exists for all the values of parameters we need. Second, we will need to show that the measures are consistent, that is, that the Chapman–Kolmogorov equations are satisfied. Finally, we need to verify that properties (1), (2) and (3) hold.

For $k = 0, 1, \ldots, n$, $j = 0, 1, \ldots, k$ let

(9)
$$
\gamma_{n,k,j} = s^j \eta^j q^{(2k-1-j)j/2} \frac{[n]_q![n-k+j-1]_q!}{[n-k]_q![n-k-1]_q![j]_q![k-j]_q!},
$$

with the conventions that:

- $\gamma_{n,k,0} = \begin{bmatrix} n \\ k \end{bmatrix}_q$, $0 \leq k \leq n$,
- $\gamma_{n,n,j} = 0$ for $0 < j \leq n$, while $\gamma_{n,n,0} = 1$,
- $\gamma_{n,k,j} = 0$ if $j < 0$, or $k < 0$, or $j > k$, or $k > n$.

Define auxiliary polynomials

(10)
$$
b_k^{(n)}(y; x, s) = \sum_{j=0}^{k} \gamma_{n,k,j}Q_{k-j}(y; x, 0, s), \qquad k = 0, 1, \ldots, n.
$$

In particular, $b_n^{(n)}(y; x, s) = Q_n(y; x, 0, s)$, $b_0^{(n)}(y; x, s) = 1$ and $b_j^{(n)}(y; x, s) = 0$ if $j < 0$ or $j > n$.

We now state the main result of this section.

THEOREM 2.1. *If polynomials $\{Q_n\}$ are defined by recurrence (6) with coefficients (7), (8) and $\{b_k^{(n)}\}$ are given by (10), then*

$$(11) \qquad Q_n(z; x, u, s) = \sum_{k=0}^{n} b_{n-k}^{(n)}(y; x, s) Q_k(z; y, u, 0)$$

*holds for $n = 1, 2, \ldots$.*

We remark that in fact a more general version of (11) holds. Let

$$\widetilde{\gamma}_{n,k,j}(t) = (-\eta)^j \frac{[n]_q! [n - k + j - 1]_q!}{[n - k]_q! [n - k - 1]_q! [j]_q! [k - j]_q!} \prod_{r=k-j}^{k-1} (t - sq^r),$$

$$\widetilde{b_k^{(n)}}(y; x, t, s) = \sum_{j=0}^{k} \widetilde{\gamma}_{n,k,j}(t) Q_{k-j}(y; x, t, s), \qquad k = 0, 1, \ldots, n.$$

Then

$$Q_n(z; x, u, s) = \sum_{k=0}^{n} \widetilde{b_{n-k}^{(n)}}(y; x, t, s) Q_k(z; y, u, t).$$

The latter formula specializes to (11) when $t = 0$. For our purposes, (11) with $y = 0$ would suffice. However, our argument relies on the three step recurrence with respect to variable $y$; this would become a mysterious trick if $y$ were set to zero, so we will prove (11).

COROLLARY 2.2.

$$(12) \qquad Q_n(z; x, u, s) = \sum_{k=1}^{n} b_{n-k}^{(n)}(0; x, s)(p_k(z; u) - p_k(x; s)), \qquad n \geq 1.$$

PROOF. Since $Q_n(z; 0, u, 0) = p_n(z; u)$, clearly (11) with $y = 0$ implies $Q_n(z; x, u, s) = \sum_{k=0}^{n} b_{n-k}^{(n)}(0; x, s) p_k(z; u)$. From (6) we see that $Q_n(x; x, s, s) = 0$, so applying the same formula again we get $\sum_{k=0}^{n} b_{n-k}^{(n)}(0; x, s) p_k(x; s) = 0$. Subtracting these two sums and noting that $p_0 = 1$ we get (12). □

The following proposition reduces Theorem 2.1 to the identities between $\{b_{n-k}^{(n)}(y; x, s)\}$.

PROPOSITION 2.3.  *Suppose a family of polynomials* $\{Q_n(y; x, t, s)\}$ *is given by recurrence (6) with arbitrary coefficients* $\mathcal{A}_n$ *and* $\mathcal{B}_n$ *such that* $\mathcal{B}_0(x, t, s) = 0$. *If* $\{b_k^{(n)}(y; x, s) : n = 0, 1, \ldots, k = 0, 1, \ldots, n\}$ *satisfy equations*

$$
\begin{aligned}
(13) \quad b_k^{(n+1)}(y; x, s) &= b_{k-1}^{(n)}(y; x, s)(\mathcal{A}_{n+1-k}(y, u, 0) - \mathcal{A}_n(x, u, s)) \\
&\quad + b_{k-2}^{(n)}(y; x, s)\mathcal{B}_{n+2-k}(y, u, 0) \\
&\quad - b_{k-2}^{(n-1)}(y; x, s)\mathcal{B}_n(x, u, s) + b_k^{(n)}(y; x, s), \\
&\hspace{4cm} k = 0, 1, 2, \ldots, n+1,
\end{aligned}
$$

*with the initial conditions* $b_0^{(0)} = 1$ *and* $b_k^{(0)} = 0$ *for* $k \geq 1$, *then (11) holds true.*

PROOF.  We proceed by induction on $n$. Obviously, (11) holds for $n = 0$ and it is also easy to verify for $n = 1$. Suppose now that $n$ is a positive integer such that (11) holds for all positive integers not greater than $n$.

By (6) and inductive assumption

$$
\begin{aligned}
Q_{n+1}&(z; x, u, s) \\
&= zQ_n(z; x, u, s) - \mathcal{A}_n(x, u, s)Q_n(z; x, u, s) - \mathcal{B}_n(x, u, s)Q_{n-1}(z; x, u, s) \\
&= z \sum_{k=0}^{n} b_{n-k}^{(n)}(y; x, s)Q_k(z; y, u, 0) \\
&\quad - \mathcal{A}_n(x, u, s) \sum_{k=0}^{n} b_{n-k}^{(n)}(y; x, s)Q_k(z; y, u, 0) \\
&\quad - \mathcal{B}_n(x, u, s) \sum_{k=0}^{n-1} b_{n-1-k}^{(n-1)}(y; x, s)Q_k(z; y, u, 0).
\end{aligned}
$$

Again by (6), the right-hand side equals

$$
\begin{aligned}
\sum_{k=0}^{n} b_{n-k}^{(n)}&(y; x, s)\{Q_{k+1}(z; y, u, 0) + \mathcal{A}_k(y, u, 0)Q_k(z; y, u, 0) \\
&\hspace{3cm} + \mathcal{B}_k(y, u, 0)Q_{k-1}(z; y, u, 0)\} \\
&\quad - \mathcal{A}_n(x, u, s) \sum_{k=0}^{n} b_{n-k}^{(n)}(y; x, s)Q_k(z; y, u, 0) \\
&\quad - \mathcal{B}_n(x, u, s) \sum_{k=0}^{n-1} b_{n-1-k}^{(n-1)}(y; x, s)Q_k(z; y, u, 0).
\end{aligned}
$$

Collecting the coefficients at $Q_j(z; y, u, 0)$ we get

$$Q_{n+1}(z; x, u, s)$$
$$= \sum_{k=0}^{n+1} Q_{n-k+1}(z; y, u, 0)\{b_k^{(n)}(y; x, s)$$
$$+ \mathcal{A}_{n-k+1}(y, u, 0)b_{k-1}^{(n)}(y; x, s)$$
$$+ \mathcal{B}_{n-k+2}(y, u, 0)b_{k-2}^{(n)}(y; x, s)$$
$$- \mathcal{A}_n(x, u, s)b_{k-1}^{(n)}(y; x, s)$$
$$- \mathcal{B}_n(x, u, s)b_{k-2}^{(n-1)}(y; x, s)\}.$$

Thus (11) follows from the assumed recurrences for $b_k^{(n+1)}$ (recall the convention $b_j^{(n+1)} = 0$ when $j < 0$ or $j > n+1$). $\square$

PROOF OF THEOREM 2.1. By Proposition 2.3, in order to prove Theorem 2.1 it suffices to verify that functions $\{b_k^{(n)}\}$ as defined by formula (10) satisfy recurrence (13). This in turn consists of expanding both sides of (13) using (10), then applying (6) to $Q_j(y; x, 0, s)$, and finally matching the coefficients at $Q_j(y; x, 0, s)$. The last task is tedious so we postpone it to the Appendix. $\square$

We now use Theorem 2.1 to prove our main result. The case $q = 0$ is covered by our argument when $0^k$ is interpreted correctly as 0 or 1; however this case is much simpler and is available separately in [12]. The case $q = -1$ is elementary and is handled separately in Section 3.2. The case $q = 1$ requires additional modifications; since an elementary self-contained construction based on birth/death processes is available in [11], in Section 3.1 we identify the transition probabilities and omit the rest of the proof.

2.1. *Proof of Theorem 1.2 for $0 < |q| < 1$.* If $\theta\eta = 0$ then $(X_t)$ is a $q$-Brownian motion when both parameters are zero, a $q$-Poisson type process when $\eta = 0$, $\theta \neq 0$, or a time inverse of a $q$-Poisson type process when $\theta = 0$, $\eta \neq 0$, see [10].

We will be therefore interested only in the case $\theta\eta \neq 0$. Passing to $(-X_t)$ if necessary, without loss of generality we may assume $\eta > 0$. Replacing our process by the process $(\sqrt{\eta/|\theta|}X_{t|\theta|/\eta})_{t \geq 0}$, without loss of generality we may further assume that $\theta = \pm\eta$.

The plan of the proof is similar to the proof for $q = 0$ as given in [12]. We first verify that $\mathcal{B}_n(x, s, t) \geq 0$ on the support of $X_s$, so that we can define the transition probabilities $\mathbf{P}_{s,t}(x, dy)$ as the orthogonality measures of $\{Q_n(y; x, t, s)\}$. We will then verify that these probability measures $\mathbf{P}_{s,t}(x, dy)$ satisfy the Chapman–Kolmogorov equation. Since $p_n(x; t) = Q_n(x; 0, t, 0)$, the second property will imply that $\{p_n(x; t)\}$ are orthogonal martingale polynomials.

Through the remainder of the proof we fix

$$(14) \qquad |q| < 1, \ \eta > 0, \text{ and } \theta \text{ such that } 1 + \eta\theta > \max\{q, 0\}.$$

2.1.1. *Support of $X_t$.* Under assumption (14) the three step recurrence (4) has positive coefficient at $p_{n-1}$, and the coefficients of the recurrence are bounded uniformly in $n$. Therefore, see [13], Theorems 2.5.4 and 2.5.5, polynomials $\{p_n(x; t)\}$ are orthogonal with respect to the unique compactly supported probability measure $\pi_t(dx)$.

We determine $\pi_t$ by setting $b = c = 0$ in [3], Theorem 2.5. To cast our recurrence (4) into the appropriate form, we need to re-parameterize our polynomials. With $t > 0$ being fixed, let $P_n(x) = p_n((\alpha x + \beta); t)/\alpha^n$, where $\alpha = \frac{2\sqrt{t}}{1-q}\sqrt{\eta\theta + 1 - q}$, $\beta = \frac{\theta + t\eta}{1-q}$. Then (4) translates into

$$xP_n(x) = P_{n+1}(x) + \tfrac{1}{2}(a+d)q^n P_n(x) + \tfrac{1}{4}(1 - adq^{n-1})(1 - q^n)P_{n-1}(x),$$

with

$$ad = \frac{\eta\theta}{\eta\theta + 1 - q}, \qquad a + d = -\frac{1}{\sqrt{t}}\frac{\theta + t\eta}{\sqrt{\eta\theta + 1 - q}}.$$

Solving this system of equations we get

$$a = -\frac{\max\{t\eta, \theta\}}{\sqrt{t}\sqrt{\eta\theta + 1 - q}},$$

$$d = -\frac{\min\{t\eta, \theta\}}{\sqrt{t}\sqrt{\eta\theta + 1 - q}}.$$

Thus $\{P_n\}$ are the Al-Salam–Chihara polynomials.

Taking into account our linear change of variable $\alpha x + \beta$, from the general theory of Askey–Wilson polynomials [3], Theorem 2.5, see also [2], we see that the support of the absolutely continuous part of $\pi_t$ is the interval $[\beta - \alpha, \beta + \alpha]$, that is,

$$(15) \qquad \left[\frac{\theta + t\eta - 2\sqrt{t}\sqrt{\eta\theta + 1 - q}}{1-q}, \frac{\theta + t\eta + 2\sqrt{t}\sqrt{\eta\theta + 1 - q}}{1-q}\right]$$

for all $t > 0$. (The absolutely continuous part exists for all $t > 0$ except when $\eta\theta + 1 - q = 0$, or when $q = 1$.)

For $0 < t < \frac{\theta^2}{\eta\theta + 1 - q}$, the discrete part of $\pi_t$ is supported on the finite number of points

$$(16) \qquad x_k = -\frac{1}{1-q}\left(\theta q^k + t\frac{\eta\theta + 1 - q}{\theta q^k} - (t\eta + \theta)\right)$$

for $k = 0, 1, \ldots$ such that $t(\eta\theta + 1 - q) < q^{2k}\theta^2$.

For $t > \frac{\eta\theta + 1 - q}{\eta^2}$, the discrete part of $\pi_t$ is supported on the finite number of points

$$(17) \qquad x_k = -\frac{1}{1-q}\left(t\eta q^k + \frac{\eta\theta + 1 - q}{\eta q^k} - (t\eta + \theta)\right),$$

for $k = 0, 1, \ldots$ such that $t\eta^2 q^{2k} > \eta\theta + 1 - q$.

These calculations lead to the following.

LEMMA 2.4. *For fixed* $0 \le t < u$ *and* $\mathcal{B}_n$ *given by* (8) *let*

$$(18) \qquad U_t = \bigcap_{n=1}^{\infty}\left\{x : \prod_{k=1}^{n} \mathcal{B}_n(x, u, t) \ge 0\right\}.$$

*Then* $\pi_t(U_t) = 1$.

(Notice that set $U_t$ does not depend on variable $u$.)

PROOF. We first verify that the absolutely continuous part of $\pi_t$ charges all of its mass to $U_t$. To see this, we show that it assigns no probability to any of the sets $V_n = \{x : 1 + \eta q^{n-1}x + [n-1]_q\eta(\theta - t\eta q^{n-1}) < 0\}$. Since $\mathcal{B}_n$ is a linear function of $x$, to verify the latter, we only need to verify that none of the endpoints of the interval (15) is in $V_n$ for all $n \ge 1$, that is,

$$1 + \eta q^{n-1}\frac{\theta + t\eta \pm 2\sqrt{t}\sqrt{\eta\theta + 1 - q}}{1 - q} + [n-1]_q\eta(\theta - t\eta q^{n-1}) \ge 0.$$

We simplify this inequality to

$$(\eta\theta + 1 - q) + t\eta^2 q^{2n-2} \ge \pm 2\eta q^{n-1}\sqrt{t}\sqrt{\eta\theta + 1 - q},$$

which holds true by the elementary inequality $A^2 + B^2 \ge \pm 2AB$.

We now verify that the discrete part of $\pi_t$ does not charge the set $U_t' = \bigcup_{n=1}^{\infty} V_n$. If $t = 0$ then $\pi_0 = \delta_0$, so $\pi_0(V_n) > 0$ implies $1 + \eta\theta[n-1]_q < 0$, contradicting (14). Suppose $t > 0$ and $x_k$ is in the support of the discrete part and is given by (16). Then $t(\eta\theta + 1 - q) < q^{2k}\theta^2$ and the condition $\mathcal{B}_{n+1}(x_k, u, t) \ge 0$ is equivalent to

$$(19) \qquad 1 + \eta q^n x_k + [n]_q\eta(\theta - t\eta q^n) \ge 0.$$

A calculation shows that this inequality is equivalent to

$$\left(1 - \frac{q^{n-k}t\eta}{\theta}\right)(1 - q + \eta\theta(1 - q^{n+k})) \ge 0.$$

Notice that $1 - q + \eta\theta(1 - q^{n+k}) = (1 - q)(1 + \eta\theta[n+k]_q) \ge 0$, as the parameters are chosen to ensure that the coefficient at the third term of recurrence (4) is nonnegative. So to show that (19) holds we only need to verify that

$$1 - \frac{q^{n-k}t\eta}{\theta} \ge 0.$$

If $q^{n-k}\theta \leq 0$ then both terms are positive and (19) holds. If $q^{n-k}\theta > 0$ then $t < q^{2k}\theta^2/(\eta\theta + 1 - q)$ implies

$$1 - \frac{q^{n-k}t\eta}{\theta} > 1 - q^{n+k}\frac{\eta\theta}{\eta\theta + 1 - q} = (1-q)\frac{1 + \eta\theta[n+k]_q}{\eta\theta + 1 - q} \geq 0,$$

so (19) follows.

Suppose now that $x_k$ is in the support of the discrete part and is given by (17). Then the condition $\mathcal{B}_{n+1}(x_k, u, t) \geq 0$ is equivalent to

$$(20) \qquad\qquad \frac{1 - q^{n-k}}{1-q}(1 - q + \eta\theta - t\eta^2 q^{n+k}) \geq 0.$$

In particular, $\mathcal{B}_{k+1}(x_k, u, t) = 0$, so we only need to show that $\mathcal{B}_{n+1}(x_k, u, t) \geq 0$ for $0 \leq n < k$.

Assuming $n < k$ we consider two cases. Suppose first that $k - n > 0$ is even or $q > 0$. Then the first factor in (20) is negative and the second factor is $(-t\eta^2 q^{2k} + 1 - q + \eta\theta) + \eta^2 t q^{2k}(1 - q^{n-k}) < 0$ as the sum of two negative terms. Suppose now that $k - n > 0$ is odd and $-1 < q < 0$. Then both factors in (20) are positive.   $\square$

### 2.1.2. *Transition probabilities.*

This part of the proof follows closely [10] and [12]. For $0 \leq s < t$ consider the six-parameter family of polynomials $\{Q_n(y; x, t, s) : n \geq -1\}$ in variable $y$ which is defined by the three step recurrence (6) with the coefficients (7), (8). (In this notation, we suppress the dependence of $Q_n$ on the parameters $\eta, \theta, q$, which are fixed.)

We define $\mathbf{P}_{s,t}(x, dy)$ as the (unique) probability measure which makes the polynomials $\{Q_n(y; x, t, s) : n \in \mathbb{N}\}$ orthogonal. This is possible whenever $1 + \eta\theta \geq \max\{q, 0\}$, $-1 < q < 1$, and $x \in U_s$ as defined by (18); the latter holds true for all $x$ from the support of the probability measure $\pi_s(dy) = \mathbf{P}_{0,s}(0, dy)$, see Lemma 2.4. Since the coefficients of the three step recurrence (6) are bounded, it is well known that measures $\{\mathbf{P}_{s,t}(x, dy)\}$ have bounded support.

The proof of Favard's theorem as given in [13], Section 2.5, shows that $\mathbf{P}_{s,t}(x, dy)$ is a weak limit of the discrete measures supported on the zeros of polynomials $\{Q_n(y; x, t, s)\}$; since these zeros depend continuously on $x$, functions $x \mapsto \mathbf{P}_{s,t}(x, A)$ are measurable.

We now follow the argument from [10], Proposition 3.2, and verify that probability measures $\{\mathbf{P}_{s,t}(x, dy)\}$ are the transition probabilities of a Markov process.

We begin by noting that since $\int Q_n(z; x, u, s)\mathbf{P}_{s,u}(x, dz) = 0$ for $n \geq 1$ and $b_0^{(n)}(y; x, s) = 1$, from (12) we recurrently get

$$(21) \qquad\qquad \int p_n(z; u)\mathbf{P}_{s,u}(x, dz) = p_n(x; s), \qquad n \geq 1.$$

PROPOSITION 2.5.  *If $0 \leq s < t$, $1 + \eta\theta \geq \max\{q, 0\}$, $-1 < q < 1$, then for a set of $x$ of $\pi_s$-measure one,*

$$(22) \qquad \mathbf{P}_{s,u}(x, \cdot) = \int_{U_t} \mathbf{P}_{t,u}(y, \cdot) \mathbf{P}_{s,t}(x, dy).$$

PROOF.  We first verify the special case $s = 0$ of (22), where we rename the variables for ease of use in the remaining proof:

$$(23) \qquad \pi_t(A) = \int_{U_s} \mathbf{P}_{s,t}(z, A) \pi_s(dz).$$

To see this, let $\nu(A) = \int_{U_s} \mathbf{P}_{s,t}(x, A) \pi_s(dx)$. Since $\mathbf{P}_{s,t}(x, \cdot)$ is a well-defined non-negative measure for all $x \in U_s$, to verify (23) we only need to show that the recurrence (4) defines polynomials $\{p_n(x; t) : n \geq 0\}$ which are orthogonal with respect to $\nu$. Since $\pi_s(U_s) = 1$, this follows from (12) in a manner similar to the general case below, so the argument is omitted.

Let $V_{s,t} = \{x \in U_s : \mathbf{P}_{s,t}(x, U_t) = 1\}$. From (23) we see that $\pi_s(V_{s,t}) = 1$. Indeed, suppose $\pi_s(V_{s,t}) < 1$ so that $\pi_s(U_s \setminus V_{s,t}) > 0$. Since $\mathbf{P}_{s,t}(x, U_t) < 1$ on $U_s \setminus V_{s,t}$, we have $1 = \pi_t(U_t) = \pi_s(V_{s,t}) + \int_{U_s \setminus V_{s,t}} \mathbf{P}_{s,t}(z, U_t) \pi_s(dz) < \pi_s(V_{s,t}) + \pi_s(U_s \setminus V_{s,t}) = 1$, a contradiction.

Let

$$\nu(A) = \int_{U_t} \mathbf{P}_{t,u}(y, A) \mathbf{P}_{s,t}(x, dy).$$

We will show that $\nu(dz) = \mathbf{P}_{s,u}(x, dz)$ for all $x \in V_{s,t}$ by verifying that the polynomials $\{Q_n(z; x, u, s)\}$ are orthogonal with respect to $\nu(dz)$. Since polynomials $\{Q_n\}$ satisfy the three step recurrence (6) with bounded coefficients, they determine a unique probability measure; since $\nu(\cdot)$ is a probability measure, to verify that it coincides with $\mathbf{P}_{s,u}(x, dz)$ it suffices to show that for $n \geq 1$ measure $\nu(dz)$ integrates $Q_n(z; x, u, s)$ to zero.

We use the fact that $\int_{\mathbb{R}} Q_n(z; y, u, t) \mathbf{P}_{t,u}(y, dz) = 0$ for $n \geq 1$. From (12) and (21) we get

$$\int_{\mathbb{R}} Q_n(z; x, u, s) \nu(dz)$$

$$= \int_{U_t} \sum_{k=1}^{n} b_{n-k}^{(n)}(0; x, s) \int_{\mathbb{R}} (p_k(z; u) - p_k(x; s)) \mathbf{P}_{t,u}(y, dz) \mathbf{P}_{s,t}(x, dy)$$

$$= \sum_{k=1}^{n} b_{n-k}^{(n)}(0; x, s) \int_{U_t} (p_k(y; t) - p_k(x; s)) \mathbf{P}_{s,t}(x, dy)$$

$$= \sum_{k=1}^{n} b_{n-k}^{(n)}(0; x, s) \int_{\mathbb{R}} (p_k(y; t) - p_k(x; s)) \mathbf{P}_{s,t}(x, dy) = 0,$$

as for $x \in V_{s,t}$ we have $\mathbf{P}_{s,t}(x, \mathbb{R} \setminus U_t) = 0$.    $\square$

The following lemma must be known but we could not find an appropriate reference.

LEMMA 2.6.    *If the family of probability measures* $\{\mathbf{P}_{s,t}(x, dy) : 0 \leq s < t\}$ *satisfies (22) on the set of $x$'s of* $\mathbf{P}_{0,s}(0, dx)$*-measure one, then there is a* $U_t$*-valued Markov process* $(X_t)$ *with* $X_0 = 0$ *and transition probabilities* $\mathbf{P}_{s,t}(x, dy)$.

CONCLUSION OF THE PROOF OF THEOREM 1.2.    By Proposition 2.5 and Lemma 2.6, there is a Markov process $(X_t)$ such that $X_0 = 0$ with transition probabilities $\mathbf{P}_{s,t}(x, dy)$.

Since $p_n(x; t) = Q_n(x; 0, t; 0)$, the construction shows that $\{p_n(x; t)\}$ are orthogonal polynomials with respect to measure $\pi_t(dx) = \mathbf{P}_{0,t}(0, dx)$. For a Markov process, the martingale polynomial property follows from (21).    $\square$

**3. Some special cases.**    In this section we complete the proof of Theorem 1.2 by analyzing two previously omitted extreme values of $q$.

3.1. *Bi-Poisson process with $q = 1$.*    In this section we read out the transition probabilities determined by (6) in the case $q = 1$; then the remaining parameters satisfy $\eta\theta \geq 0$. This case is not covered by the construction in Section 2.1 as the coefficients of recurrence (4) are unbounded. After we identify the transition probabilities from (6), we will omit the construction, as an even more elementary construction based on birth/death processes is presented in [11].

PROPOSITION 3.1.    *Suppose* $\eta, \theta > 0$ *and* $q = 1$. *If* $(Y_t)$ *is a Markov process with transition probabilities such that polynomials defined by the three step recurrence (6) with coefficients (7), (8) are orthogonal, then*

$$Y_t = \begin{cases} (\theta - \eta t)Z_t - \dfrac{t}{\theta}, & 0 \leq t < \theta/\eta, \\[2mm] Z_{\theta/\eta} - \dfrac{1}{\eta}, & t = \theta/\eta, \\[2mm] (\eta t - \theta)Z_t - \dfrac{1}{\eta}, & t > \theta/\eta, \end{cases}$$

*where* $(Z_t)$ *is Markov with the transition probabilities determined as follows.*

$$Z_t - Z_s | Z_s \overset{\mathcal{D}}{=} \mathrm{NB}\left(\frac{1}{\theta\eta} + Z_s, \frac{\theta - \eta t}{\theta - \eta s}\right), \qquad 0 \leq s < t < \theta/\eta,$$

$$Z_{\theta/\eta} | Z_s \overset{\mathcal{D}}{=} \mathrm{Gam}\left(\frac{1}{\theta\eta} + Z_s, \theta - \eta s\right), \qquad 0 \leq s < \theta/\eta,$$

Table 1
*Distributions of the bi-Poisson process with $q = 1$*

| Name | Parameters | Distribution | $E(e^{uZ})$ | Notation |
|------|-----------|-------------|-------------|----------|
| Poisson | $\lambda > 0$ | $e^{-\lambda}\lambda^k/k!, k = 0, 1, 2, \ldots$ | $\exp(\lambda(e^u - 1))$ | $\mathrm{Poiss}(\lambda)$ |
| Gamma | $r > 0, \sigma > 0$ | $f(z) = \frac{x^{r-1}}{\sigma^r \Gamma(r)} e^{-x/\sigma}$ | $(1 - \sigma u)^{-r}$ | $\mathrm{Gam}(r, \sigma)$ |
| Negative Binomial | $r > 0, 0 < p < 1$ | $\frac{\Gamma(k+r)}{\Gamma(r)k!} p^r(1-p)^k, k = 0, 1, \ldots$ | $\frac{p^r}{(1-(1-p)e^u)^r}$ | $\mathrm{NB}(r, p)$ |
| Binomial | $n \geq 0, 0 \leq p \leq 1$ | $\binom{n}{k} p^k (1-p)^{n-k}, k = 0, 1, \ldots, n.$ | $(1 - p + pe^u)^n$ | $\mathrm{Bin}(n, p)$ |

$$Z_t | Z_{\theta/\eta} \overset{\mathcal{D}}{=} \mathrm{Poiss}\left( \frac{Z_{\theta/\eta}}{\eta t - \theta} \right), \qquad t > \theta/\eta,$$

$$Z_t | Z_s \overset{\mathcal{D}}{=} \mathrm{Bin}\left( Z_s, \frac{\eta s - \theta}{\eta t - \theta} \right), \qquad \theta/\eta < s < t.$$

*Here we use the fairly standard notation for the distributions summarized in Table 1.*

The following result can be traced back to [15]; we state it in the parametrization used in [10].

LEMMA 3.2. *Fix parameters $\tilde{t} > 0$, $\tilde{\theta}$, $\tilde{\tau} \in \mathbb{R}$. Suppose polynomials $\{\tilde{p}_n(y) : n \geq 0\}$ are given by the three step recurrence*

$$(24) \qquad y\tilde{p}_n(y) = \tilde{p}_{n+1}(y) + \tilde{\theta}n\tilde{p}_n(y) + (\tilde{t} + \tilde{\tau}(n-1))n\tilde{p}_{n-1}(y)$$

*for $n = 1, 2, \ldots$ with $\tilde{p}_{-1} = 0$ and $\tilde{p}_0 = 1$.*

*Let $Y$ be a random variable such that $\{\tilde{p}_n(Y)\}$ are orthogonal.*

(i) *If $\tilde{\tau} = 0$, $\tilde{\theta} \neq 0$, then*

$$Y = \tilde{\theta}Z - \tilde{t}/\tilde{\theta},$$

*where $Z$ is Poisson with parameter $\lambda = \tilde{t}/\tilde{\theta}^2$.*

(ii) *If $\tilde{\theta}^2 = 4\tilde{\tau} > 0$ then*

$$Y = \mathrm{sgn}(\tilde{\theta})Z - 2\tilde{t}/\tilde{\theta},$$

*where $Z$ is Gamma with shape parameter $r = \tilde{t}/\tilde{\tau}$ and scale parameter $\sigma = |\tilde{\theta}|/2$.*

(iii) *If $\tilde{\theta}^2 > 4\tilde{\tau} > 0$ then*

$$Y = \mathrm{sgn}(\tilde{\theta})\sqrt{\tilde{\theta}^2 - 4\tilde{\tau}}Z - \mathrm{sgn}(\tilde{\theta})\frac{|\tilde{\theta}| - \sqrt{\tilde{\theta}^2 - 4\tilde{\tau}}}{2\tilde{\tau}}\tilde{t},$$

*where $Z$ is negative binomial with parameters $r = \tilde{t}/\tilde{\tau}$, $p = \frac{2\sqrt{\tilde{\theta}^2 - 4\tilde{\tau}}}{|\tilde{\theta}| + \sqrt{\tilde{\theta}^2 - 4\tilde{\tau}}}$.*

(iv) *If $\tilde\tau < 0$ and $n = -\tilde t/\tilde\tau$ is a natural number then*

$$Y = \sqrt{\tilde\theta^2 - 4\tilde\tau}\, Z + \frac{\tilde t}{2\tilde\tau}(\sqrt{\tilde\theta^2 - 4\tilde\tau} - \tilde\theta),$$

*where $Z$ is binomial $Bin(n, p)$ with $p = \frac{1}{2}(1 - \frac{\tilde\theta}{\sqrt{\tilde\theta^2 - 4\tilde\tau}})$.*

NOTE.  If (24) holds with $\tilde t = 0$ then $Y = 0$.

PROOF.  The type of the law of $Y$ can be read out from [16], Table 1 and formula (8.3). Perhaps the simplest way to identify the parameters of the basic law of $Z$ as well as the linear transformation $Y = a(Z - \mathbb{E}(Z))$, is to compare the first four moments. We have

$$\mathbb{E}(Y) = 0, \qquad \mathbb{E}(Y^2) = \tilde t, \qquad \mathbb{E}(Y^3) = \tilde\theta\tilde t, \qquad \mathbb{E}(Y^4) = 3\tilde t^2 + (\tilde\theta^2 + 2\tilde\tau)\tilde t.$$

Solving the system of equations $a^j \mathbb{E}(Z - \mathbb{E}Z)^j = \mathbb{E}(Y^j)$, $j = 2, 3, 4$, we get the results as stated.

Except for (iv), the type of the distribution can also be read out from [10], Theorem 4.2. However caution is needed as some of the explicit functions are typeset incorrectly in the published version (the formulas are corrected in the arXiv preprint).  □

PROOF OF PROPOSITION 3.1.  Fixed $0 \le s < t$. Recurrence (6) for polynomials $Q_n(z + x; x, t, s)$ is equivalent to (24) with the following identification of parameters:

$$\tilde\theta = t\eta + \theta - 2s\eta, \qquad \tilde\tau = \eta(\theta - s\eta)(t - s), \qquad \tilde t = (1 + \eta x)(t - s).$$

A calculation shows that $\tilde\theta^2 - 4\tilde\tau = (t\eta - \theta)^2$. Thus if $Y_s > -1/\eta$ then $\tilde t > 0$ and one of the cases listed in Lemma 3.2 must hold.

We now use (6) with $x = s = 0$ to verify that for $t > 0$ we have $Y_t \ge -1/\eta$. For $t \neq \theta/\eta$ by Lemma 3.2 we have

$$Y_t = |\theta - \eta t| Z_t - \frac{\min\{\theta, \eta t\}}{\theta\eta},$$

where $Z_t \sim NB(\frac{1}{\eta\theta}, \frac{|\theta - \eta t|}{\max\{\theta, \eta t\}})$. For $t = \theta/\eta$ we have $Y_{\theta/\eta} \stackrel{d}{=} Z - \frac{1}{\eta}$, where $Z \sim$ Gamma$(\frac{1}{\theta\eta}, \theta)$. Therefore, $Y_s \ge -1/\eta$ and the inequality is strict except when $s > \theta/\eta$.

We now use Lemma 3.2 to identify the conditional distribution of $Y_t$ given $Y_s$.

(i) Suppose $0 \le s < t < \theta/\eta$. Then $\tilde\theta > \theta - t\eta > 0$ and $\tilde\tau > 0$.

$$Y_t - Y_s | Y_s \stackrel{\mathcal{D}}{=} (\theta - t\eta)Z - (t - s)\frac{1 + \eta Y_s}{\theta - s\eta},$$

where $Z|Y_s$ is negative binomial with parameters $r = \frac{Y_s + 1/\eta}{\theta - s\eta}$ and $p = \frac{\theta - t\eta}{\theta - s\eta}$. Writing $Y_s = (\theta - \eta s)Z_s - s/\theta$ we get $r = \frac{1}{\eta\theta} + Z_s$ and

$$Y_t - Y_s|Y_s \overset{\mathcal{D}}{=} (\theta - t\eta)Z - (t - s)/\theta - (t - s)Z_s\eta.$$

Letting $Z = Z_t - Z_s$ we get

$$Y_t|Y_s \overset{\mathcal{D}}{=} (\theta - t\eta)Z_t - t/\theta.$$

(ii) Suppose $0 \le s < t = \theta/\eta$. Then $\tilde{\theta}^2 - 4\tilde{\tau} = 0$ and $\tilde{\tau} > 0$. Thus $Y_t - Y_s|Y_s \overset{\mathcal{D}}{=} Z_t - Y_s - 1/\eta$, so $Y_t|Y_s \overset{\mathcal{D}}{=} Z_t - 1/\eta$, where $Z_t|Y_s$ is Gamma with $r = \frac{Y_s + 1/\eta}{\theta/\eta - s}$ and $\sigma = \theta - s\eta$. Writing $Y_s = (\theta - \eta s)Z_s - s/\theta$ we get $r = \frac{1}{\eta\theta} + Z_s$.

(iii) Suppose $s = \theta/\eta < t$. Then $\tilde{\tau} = 0$, $\tilde{\theta} = t\eta - \theta > 0$. Thus $Y_t - Y_s|Y_s \overset{\mathcal{D}}{=} (t\eta - \theta)Z_t - (1/\eta + Y_s)$, that is, $Y_t|Y_s \overset{\mathcal{D}}{=} (t\eta - \theta)Z_t - 1/\eta$, where $Z_t|Y_s$ is Poiss($\lambda$) with $\lambda = \frac{Y_s + 1/\eta}{t\eta - \theta} = \frac{Z_s}{t\eta - \theta}$.

(iv) Suppose $\theta/\eta < s < t$. Then $\tilde{\tau} < 0$. Therefore, there must be an integer $n = Z_s$ such that $\tilde{t} = -\tilde{\tau}n$, that is, $Y_s = (s\eta - \theta)Z_s - 1/\eta$. Then $Y_t - Y_s|Y_s \overset{\mathcal{D}}{=} (t\eta - \theta)Z_t - (Y_s + 1/\eta) = (t\eta - \theta)Z_t - (s\eta - \theta)Z_s$, where $Z_t|Z_s$ is binomial with $p = \frac{s\eta - \theta}{t\eta - \theta}$ and $n = Z_s$. Thus

$$Y_t|Z_s \overset{\mathcal{D}}{=} (t\eta - \theta)Z_t - 1/\eta.$$

(v) Suppose $0 \le s < \theta/\eta < t$. Then $\tilde{\theta} = (t - s)\eta + \theta - s\eta > 0$ and $\tilde{\tau} > 0$. Therefore $Y_t - Y_s|Y_s \overset{\mathcal{D}}{=} (t\eta - \theta)Z - \frac{1 + \eta Y_s}{\eta}$, where $Z|Y_s$ is negative binomial with $r = \frac{Y_s + 1/\eta}{\theta - s\eta}$ and $p = \frac{t\eta - \theta}{\eta(t - s)}$. Thus

$$Y_t|Y_s \overset{\mathcal{D}}{=} (t\eta - \theta)Z - 1/\eta.$$

Writing $Y_s = (\theta - \eta s)Z_s - s/\theta$, $Z = Z_t$, we get $r = \frac{1}{\eta\theta} + Z_s$ and

$$Y_t|Y_s \overset{\mathcal{D}}{=} (t\eta - \theta)Z_t - 1/\eta.$$

This shows that $(Z_t)$ is a Markov process with the transition probabilities as stated. $\square$

3.2. *Bi-poisson process with $q = -1$.* In this section we read out the transition probabilities determined by (6) in the case $q = -1$; then the remaining parameters satisfy $1 + \eta\theta \ge 0$. This case is not covered by the construction in Section 2.1, as the reasoning based on [3] fails.

We first observe that the polynomials $p_n(x;t) = 0$ on the support of $\pi_t$ for $n > 1$, $t > 0$. Indeed, since $[2]_q = 0$, from (4) we see that $p_n(x;t)$ has $p_2(x;t)$ as a factor for $n \geq 2$. Thus the univariate distributions of $(X_t)$ are supported at two points $a_\pm(t)$ which are the roots of $p_2(x; t)$. The roots are

$$a_\pm(t) = \frac{t\eta + \theta \pm \sqrt{4t + (t\eta + \theta)^2}}{2}.$$

From the requirement that $p_1(x; t) = x$ is orthogonal to $p_2$ (i.e., that the mean of $X_t$ is zero), we compute the corresponding probabilities

$$p_\pm(t) := \Pr(X_t = a_\pm(t)) = \frac{1}{2} \mp \frac{t\eta + \theta}{2\sqrt{4t + (t\eta + \theta)^2}}.$$

A calculation confirms that $\mathbb{E}(X_t^2) = a_+^2 p_+ + a_-^2 p_- = t$.

Similarly, see (6), the transition probabilities that make polynomials $\{Q_n\}$ orthogonal are supported on the roots of

$$Q_2(y; x, t, s) = y^2 + y((s\eta - x)(1 + q) - t\eta - \theta) + s - t + qx^2 + x\theta - qsx\eta.$$

A calculation shows that the supports of the one-dimensional distributions are consistent:

$$Q_2(a_\pm(t); a_\pm(s), t, s) = 0.$$

We now determine the transition probabilities

$$p_{\pm\pm}(s, t) = \Pr(X_t = a_\pm(t) | X_s = a_\pm(s))$$

from the requirement that $\mathbb{E}(X_t | X_s) = X_s$. A calculation gives the transition matrix

$$P_{s,t} = \begin{bmatrix} p_{--}(s,t) & p_{-+}(s,t) \\ p_{+-}(s,t) & p_{++}(s,t) \end{bmatrix}$$

with

$$p_{++}(s,t) = \frac{s\eta - t\eta + \sqrt{4s + (s\eta + \theta)^2} + \sqrt{4t + (t\eta + \theta)^2}}{2\sqrt{4t + (t\eta + \theta)^2}},$$

$$p_{-+}(s,t) = \frac{s\eta - t\eta - \sqrt{4s + (s\eta + \theta)^2} + \sqrt{4t + (t\eta + \theta)^2}}{2\sqrt{4t + (t\eta + \theta)^2}},$$

$$p_{+-}(s,t) = \frac{-s\eta + t\eta - \sqrt{4s + (s\eta + \theta)^2} + \sqrt{4t + (t\eta + \theta)^2}}{2\sqrt{4t + (t\eta + \theta)^2}},$$

$$p_{--}(s,t) = \frac{-s\eta + t\eta + \sqrt{4s + (s\eta + \theta)^2} + \sqrt{4t + (t\eta + \theta)^2}}{2\sqrt{4t + (t\eta + \theta)^2}}.$$

Using these formulas, a calculation verifies the Chapman–Kolmogorov equation

$$P_{s,t}P_{t,u} = P_{s,u}.$$

Thus the two-valued Markov process $X_t \in \{a_+(t), a_-(t)\}$ is well defined, and has mean zero and variance $t$.

To verify (2) and (3), we compute the conditional probabilities: for $\alpha$, $\beta$, $\gamma = \pm$,

$$\Pr(X_t = a_\beta(t)|X_s = a_\alpha(s), X_u = a_\gamma(u)) = \frac{p_{\alpha\beta}(s,t)p_{\beta\gamma}(t,u)}{p_{\alpha\gamma}(s,u)}.$$

A long but straightforward calculation now verifies (2) and (3).

From $\mathbb{E}[X_t|X_s] = X_s$, we deduce $\mathbb{E}(X_s X_s) = s$ for $s < t$ so (1) holds.

## APPENDIX: COMPUTATIONS NEEDED FOR THE PROOF OF THEOREM 2.1

Our aim is to show (13) for $b_k^{(n)}$ defined in (10) and (9).

We first notice $b_0^{(n)} = \gamma_{n,0,0} Q_0(y; x, 0, s) = \gamma_{n,0,0} = 1$, so (13) holds for $k = 0$.

It remains therefore to verify (13) for $k \geq 1$. We note that (13) takes a slightly different form for the boundary values $k = 1, n, n+1$ however due to our convention these cases are covered by (13).

Since the right-hand side of (13) is an affine function of $u$ we split it into two equalities, which will be proved separately: the first for the coefficient of $u$

$$
\begin{aligned}
&[n+2-k]_q \{1 + \eta q^{n+1-k} y + [n+1-k]_q \eta\theta\} b_{k-2}^{(n)}(y; x, s) \\
&\quad = \eta[k-1]_q q^{n+1-k} b_{k-1}^{(n)}(y; x, s) \\
&\qquad + [n]_q \{1 + \eta q^{n-1} x + [n-1]_q \eta(\theta - \eta q^{n-1} s)\} b_{k-2}^{(n-1)}(y; x, s),
\end{aligned}
\tag{25}
$$

and the second for the intercept

$$
\begin{aligned}
&b_k^{(n+1)}(y; x, s) \\
&\quad = \{q^{n+1-k} y + [n+1-k]_q \theta \\
&\qquad - q^n x - [n]_q(\theta - \eta(1+q)q^{n-1}s)\} b_{k-1}^{(n)}(y; x, s) \\
&\qquad + [n]_q q^{n-1} s\{1 + \eta q^{n-1} x + [n-1]_q \eta(\theta - \eta q^{n-1} s)\} b_{k-2}^{(n-1)}(y; x, s) \\
&\qquad + b_k^{(n)}(y; x, s).
\end{aligned}
\tag{26}
$$

*First step.* We start with proving (25). We expand $b_{k-1}^{(n)}(y; x, s)$, $b_{k-2}^{(n)}(y; x, s)$ and $b_{k-2}^{(n-1)}(y; x, s)$ in (25) using (10) and then after applying the three terms recurrence formula (6) we arrive at the following equality which is equivalent to (25):

$$(27) \qquad \sum_{j=0}^{k-1} C_{n,k,j} Q_j(y; x, 0, s) = 0,$$

where

$$
\begin{aligned}
C_{n,k,j} = {} & -\eta[k-1]_q q^{n+1-k} \gamma_{n,k-1,k-1-j} \\
& + \eta[n+2-k]_q q^{n+1-k} \{ \gamma_{n,k-2,k-1-j} + \mathcal{A}_j(x,0,s)\gamma_{n,k-2,k-2-j} \\
& \hspace{6em} + \mathcal{B}_{j+1}(x,0,s)\gamma_{n,k-2,k-3-j} \} \\
& + [n+2-k]_q \{ 1 + [n+1-k]_q \eta\theta \} \gamma_{n,k-2,k-2-j} \\
& - [n]_q \{ 1 + \eta q^{n-1} x + [n-1]_q \eta(\theta - \eta q^{n-1} s) \} \gamma_{n-1,k-2,k-2-j}.
\end{aligned}
$$

We will show that $C_{n,k,j} = 0$. To this end we note that due to (7) and (8) the expression for $C_{n,k,j}$ can be written as

$$C_{n,k,j} = \eta\theta[n+2-k]_q \mathbf{a} + q^{n+1-k+j}\eta x[n+2-k]_q \mathbf{b} + q^{n+1-k}\eta\mathbf{c} + \mathbf{d},$$

where

$$
\begin{aligned}
\mathbf{a} = {} & [j]_q q^{n+1-k} \gamma_{n,k-2,k-2-j} - [j]_q[j+1]_q q^j \eta s \gamma_{n,k-2,k-3-j} \\
& + [n+1-k]_q \gamma_{n,k-2,k-2-j} - \frac{[n]_q[n-1]_q}{[n+2-k]_q} \gamma_{n-1,k-2,k-2-j},
\end{aligned}
$$

$$
\begin{aligned}
\mathbf{b} = {} & \gamma_{n,k-2,k-2-j} - [j+1]_q q^j \eta s \gamma_{n,k-2,k-3-j} \\
& - \frac{[n]_q}{[n+2-k]_q} q^{k-2-j} \gamma_{n-1,k-2,k-2-j},
\end{aligned}
$$

$$
\begin{aligned}
\mathbf{c} = {} & -[k-1]_q \gamma_{n,k-1,k-1-j} + [n]_q[n-1]_q q^{k-2} \eta s \gamma_{n-1,k-2,k-2-j} \\
& + [n+2-k]_q \{ \gamma_{n,k-2,k-1-j} - (1+q)q^{j-1}\eta s \gamma_{n,k-2,k-2-j} \\
& \hspace{6em} + [j]_q[j+1]_q q^{2j}(\eta s)^2 \gamma_{n,k-2,k-3-j} \},
\end{aligned}
$$

$$
\begin{aligned}
\mathbf{d} = {} & -[n+2-k]_q[j+1]_q q^{n+1-k+j}\eta s \gamma_{n,k-2,k-3-j} + [n+2-k]_q \gamma_{n,k-2,k-2-j} \\
& - [n]_q \gamma_{n-1,k-2,k-2-j}.
\end{aligned}
$$

We will show that $\mathbf{a} = \mathbf{b} = \mathbf{c} = \mathbf{d} = 0$.

$\mathbf{a} = 0$: From (9) we get

$$(28) \qquad [j+1]_q q^j \eta s \gamma_{n,k-2,k-3-j} = \frac{[k-2-j]_q}{[n-1-j]_q} \gamma_{n,k-2,k-2-j}$$

and

$$(29) \qquad [n]_q \gamma_{n-1,k-2,k-2-j} = \frac{[n+2-k]_q[n+1-k]_q}{[n-1-j]_q} \gamma_{n,k-2,k-2-j}.$$

Consequently

$$\mathbf{a} = \frac{[n+2-k]_q \gamma_{n,k-2,k-2-j}}{[n-1-j]_q} \{[j]_q q^{n+1-k}([n-1-j]_q - [k-2-j]_q)$$
$$- [n+1-k]_q([n-1]_q - [n-1-j]_q)\}.$$

Using twice the $q$-subtraction formula

$$(30) \qquad\qquad [m]_q - [l]_q = q^l[m-l]_q, \qquad l \le m,$$

to $[n-1-j]_q - [k-2-j]_q$, and to $[n-1]_q - [n-1-j]_q$ we get $\mathbf{a} = 0$.

$\mathbf{b} = 0$: To calculate $\mathbf{b}$ we first use (28) and (29), getting

$$\mathbf{b} = \frac{\gamma_{n,k-2,k-2-j}}{[n-1-j]_q} \{[n-1-j]_q - [k-2-j]_q - q^{k-2-j}[n+1-k]_q\}.$$

Now we use (30) to $[n-1-j]_q - [k-2-j]_q$ which gives $\mathbf{b} = 0$.

$\mathbf{c} = 0$: For $\mathbf{c}$ we need two more formulas: equation (28) with $j-1$ gives

$$(31) \qquad\qquad \gamma_{n,k-2,k-1-j} = \frac{[n-j]_q[j]_q}{[k-j-1]_q} q^{j-1} \eta s \gamma_{n,k-2,k-2-j},$$

and from (9) we get

$$(32) \qquad \gamma_{n,k-1,k-1-j} = \frac{[n+2-k]_q[n+1-k]_q}{[k-j-1]_q} q^{k-2} \eta s \gamma_{n,k-2,k-2-j}.$$

Inserting (28), (31) and (32) into the definition of $\mathbf{c}$ we get

$$\mathbf{c} = \frac{\eta s q^{n-k+1}[n+2-k]_q \gamma_{n,k-2,k-2-j}}{[n-1-j]_q[k-1-j]_q}$$
$$\times \{-[k-1]_q[n+1-k]_q[n-1-j]_q q^{k-1-j}$$
$$+ [n-j]_q[j]_q[n-1-j]_q$$
$$- [j]_q[n-1-j]_q[k-1-j]_q - [j]_q[n-1-j]_q[k-1-j]_q q$$
$$+ [j]_q[k-2-j]_q[k-1-j]_q q$$
$$+ [n-1]_q[n+1-k]_q[k-1-j]_q q^{k-1-j}\}.$$

Using (30) for $[n-j]_q - [k-1-j]_q$ we combine the second and the third expression in the brackets above, and applying (30) to $[n-1-j]_q - [k-2-j]_q$ we combine the fourth and the fifth expression in the brackets. Thus the expression in the bracket takes the form

$$q^{k-1-j}\{-[k-1]_q[n+1-k]_q[n-1-j]_q+[j]_q[n-1-j]_q[n+1-k]_q$$
$$-[j]_q[k-1-j]_q[n+1-k]_q+[n-1]_q[n+1-k]_q[k-1-j]_q\}.$$

Now we use again (30) to $[k-1]_q - [j]_q$ to combine the first two expressions, and to $[n-1]_q - [j]_q$ for the last two expression to conclude that $\mathbf{c} = 0$.

$\mathbf{d} = 0$: First apply (28) and (29) to get

$$\mathbf{d} = \frac{[n+2-k]_q\gamma_{n,k-2,k-2-j}}{[n-1-j]_q}\{[n-1-j]_q-[n+1-k]_q-q^{n+1-k}[k-2-j]_q\}.$$

Then by (30) applied to $[n-1-j]_q - [n+1-k]_q$ we see that $\mathbf{d} = 0$.

*Second step.* Next, we prove (26). We expand $b_k^{(n+1)}(y;x,s)$, $b_{k-1}^{(n)}(y;x,s)$, $b_{k-2}^{(n-1)}(y;x,s)$ and $b_k^{(n)}(y;x,s)$ in (26) using (10) and then after applying the three terms recurrence formula (6) we arrive at the following equality which is equivalent to (26):

$$(33) \qquad \sum_{j=0}^{k-1} D_{n,k,j} Q_j(y;x,0,s) = 0,$$

where

$$D_{n,k,j} = \gamma_{n+1,k,k-j} - q^{n+1-k}\gamma_{n,k-1,k-j} - q^{n+1-k}\mathcal{A}_j(x,0,s)\gamma_{n,k-1,k-1-j}$$
$$- q^{n+1-k}\mathcal{B}_{j+1}(x,0,s)\gamma_{n,k-1,k-2-j}$$
$$- \{[n+1-k]_q\theta - q^n x - [n]_q(\theta - (1+q)q^{n-1}\eta s)\}\gamma_{n,k-1,k-1-j}$$
$$- [n]_q q^{n-1}s\{1+\eta q^{n-1}x+[n-1]_q\eta(\theta - q^{n-1}\eta s)\}\gamma_{n-1,k-2,k-2-j}$$
$$- \gamma_{n,k,k-j}.$$

We will show that $D_{n,k,j} = 0$. To this end we note that due to (7) and (8) the expression for $D_{n,k,j}$ can be written as

$$D_{n,k,j} = q^{n+1-k+j}x\mathbf{A} + q^{n+1-k}\theta\mathbf{B} + q^{n+1-k+j}s\mathbf{C} + \mathbf{D},$$

where

$$\mathbf{A} = -\gamma_{n,k-1,k-1-j} + [j+1]_q q^j \eta s\gamma_{n,k-1,k-2-j}$$
$$+ q^{k-1-j}\gamma_{n,k-1,k-1-j} - [n]_q q^{n-3+k-j}\eta s\gamma_{n-1,k-2,k-2-j},$$
$$\mathbf{B} = -[j]_q\gamma_{n,k-1,k-1-j} + [j+1]_q[j]_q q^j \eta s\gamma_{n,k-1,k-2-j}$$
$$+ [k-1]_q\gamma_{n,k-1,k-1-j} - [n]_q[n-1]_q q^{k-2}\eta s\gamma_{n-1,k-2,k-2-j},$$
$$\mathbf{C} = [j+1]_q\gamma_{n,k-1,k-2-j} - [n]_q q^{k-2-j}\gamma_{n-1,k-2,k-2-j},$$

$$\mathbf{D} = \gamma_{n+1,k,k-j} - q^{n+1-k}\gamma_{n,k-1,k-j} + [j]_q(1+q)q^{n-k+j}\eta s\gamma_{n,k-1,k-1-j}$$
$$- [j+1]_q[j]_q q^{n+1-k+2j}(\eta s)^2\gamma_{n,k-1,k-2-j}$$
$$- [n]_q(1+q)q^{n-1}\eta s\gamma_{n,k-1,k-1-j}$$
$$+ [n]_q[n-1]_q q^{2n-2}(\eta s)^2\gamma_{n-1,k-2,k-2-j} - \gamma_{n,k,k-j}.$$

We will show that $\mathbf{A} = \mathbf{B} = \mathbf{C} = \mathbf{D} = 0$.

$\mathbf{A} = 0$: Equation (31) with $k$ replaced by $k+1$ gives

$$(34) \qquad [j+1]_q q^j\eta s\gamma_{n,k-1,k-2-j} = \frac{[k-1-j]_q}{[n-1-j]_q}\gamma_{n,k-1,k-1-j}.$$

Moreover by (9) we have

$$(35) \qquad [n]_q q^{k-2}\eta s\gamma_{n-1,k-2,k-2-j} = \frac{[k-1-j]_q}{[n-1-j]_q}\gamma_{n,k-1,k-1-j}.$$

Consequently, by (34) and (35), we get

$$\mathbf{A} = \frac{\gamma_{n,k-1,k-1-j}}{[n-1-j]_q}\{-[n-1-j]_q + [k-1-j]_q$$
$$+ q^{k-1-j}[n-1-j]_q - q^{n-1-j}[k-1-j]_q\}.$$

Applying now (30) to the first two terms in the brackets we get

$$-([n-1-j]_q - [k-1-j]_q) = -q^{k-1-j}[n-k]_q.$$

Combining this with the third term we get by (30) again

$$q^{k-1-j}([n-1-j]_q - [n-k]_q) = q^{n-1-j}[k-1-j]_q.$$

Thus together with the fourth term we arrive at $\mathbf{A} = 0$.

$\mathbf{B} = 0$: By (34) and (35) we get

$$\mathbf{B} = \frac{\gamma_{n,k-1,k-1-j}}{[n-1-j]_q}\{-[j]_q[n-1-j]_q + [j]_q[k-1-j]_q$$
$$+ [k-1]_q[n-1-j]_q - [n-1]_q[k-1-j]_q\}.$$

Using (30) we have $[n-1-j]_q - [k-1-j]_q = q^{k-1-j}[n-k]_q$, which is applied to the first two terms in curly braces above and also to the third term. Thus the expression in curly braces takes the form

$$-[j]_q[n-k]_q q^{k-1-j} + [k-1]_q[n-k]_q q^{k-1-j}$$
$$+ [k-1]_q[k-1-j]_q - [n-1]_q[k-1-j]_q.$$

Now combining the first two terms by (30) we get $[n-k]_q[k-1-j]_q q^{k-1}$ and combining the last two terms, again by (30) we obtain $-[k-1-j]_q[n-k]_q q^{k-1}$, and thus $\mathbf{B} = 0$.

$\mathbf{C} = 0$: By (9) it follows that

$$[n]_q q^{k-2} \gamma_{n-1,k-2,k-2-j} = [j+1]_q q^j \gamma_{n,k-1,k-2-j}$$

which implies immediately $\mathbf{C} = 0$.

$\mathbf{D} = 0$: By (9), see (35),

$$(36) \qquad \gamma_{n+1,k,k-j} = \frac{[n+1]_q [n-j]_q}{[k-j]_q} q^{k-j} \eta s \gamma_{n,k-1,k-1-j}.$$

From (28) and (29) with $k$ replaced by $k+1$, we get, respectively

$$(37) \qquad \gamma_{n,k-1,k-j} = \frac{[n-j]_q [j]_q}{[k-j]_q} q^{j-1} \eta s \gamma_{n,k-1,k-1-j},$$

$$(38) \qquad \gamma_{n,k,k-j} = \frac{[n+1-k]_q [n-k]_q}{[k-j]_q} q^{k-1} \eta s \gamma_{n,k-1,k-1-j}.$$

Plugging (34)–(38) to the definition of $\mathbf{D}$ we can express each term through $\gamma_{n,k-1,k-1-j}$ and thus

$$\begin{aligned}
\mathbf{D} = \frac{\eta s \gamma_{n,k-1,k-1-j}}{[k-j]_q [n-1-j]_q} \{ & [n+1]_q [n-j]_q [n-1-j]_q q^{k-1} \\
& - [n-j]_q [j]_q [n-1-j]_q q^{n-k+j} \\
& + [j]_q [k-j]_q [n-1-j]_q (1+q) q^{n-k+j} \\
& - [j]_q [k-j]_q [k-1-j]_q q^{n+1-k+j} \\
& - [n]_q [n-1-j]_q [k-j]_q (1+q) q^{n-1} \\
& + [n-1]_q [k-j]_q [k-1-j]_q q^{2n-k} \\
& - [n+1-k]_q [n-k]_q [n-1-j]_q q^{k-1} \}.
\end{aligned}$$

Now we combine the fourth and the sixth term in the brackets using (30) three times:

$$\begin{aligned}
[k-j]_q & q^{n+1-k+j} \{ [n-1]_q [k-1-j]_q q^{n-1-j} - [k-1-j]_q [j]_q \} \\
= [k-j]_q & q^{n+1-k+j} \{ ([j]_q + q^j [n-1-j]_q) [k-1-j]_q q^{n-1-j} \\
& - ([n-1-j]_q - q^{k-1-j}[n-k]_q) [j]_q \} \\
= [k-j]_q & \{ [n-1-j]_q [k-1-j]_q q^{2n+j-k} - [n-1-j]_q [j]_q q^{n+1-k+j} \\
& + [j]_q q^n ([n-k]_q + q^{n-k}[k-1-j]_q) \} \\
= [k-j]_q & [n-1-j]_q ([k-1-j]_q q^{2n-k+j} + [j]_q (q^n - q^{n+1-k+j})).
\end{aligned}$$

After inserting this in the last expression for $\mathbf{D}$ we note that the factor $[n-1-j]_q$ appears in every summand in the brackets. Canceling it we get

$$\mathbf{D} = \frac{\eta s \gamma_{n,k-1,k-1-j}}{[k-j]_q} \{ [n+1]_q[n-j]_q q^{k-1} - [n-j]_q[j]_q q^{n-k+j}$$
$$+ [j]_q[k-j]_q q^{n-k+j} + [k-j]_q[k-1-j]_q q^{2n-k+j}$$
$$+ [k-j]_q[j]_q q^n$$
$$- [n]_q[k-j]_q q^{n-1} - [n]_q[k-j]_q q^n$$
$$- [n+1-k]_q[n-k]_q q^{k-1} \}.$$

Adding the second and the third term in the brackets above we get by (30)

$$-[j]_q q^{n-k+j}([n-j]_q - [k-j]_q) = -[j]_q[n-k]_q q^n$$

and adding to this sum the eighth term we obtain

$$-[n-k]_q q^{k-1}([n+1-k]_q + [j]_q q^{n+1-k}) = -[n-k]_q[n+1-k+j]_q q^{k-1}.$$

Now we combine this sum with the first term:

$$q^{k-1}\{[n+1]_q[n-j]_q - [n-k]_q[n+1-k+j]_q\}$$
$$= q^{k-1}\{([n+1-k+j]_q + q^{n+1-k+j}[k-j]_q)[n-j]_q$$
$$- ([n-j]_q - q^{n-k}[k-j]_q)[n+1+k-j]_q\}$$
$$= [k-j]_q([n-j]_q q^{n+j} + [n+1-k+j]_q q^{n-1}).$$

Thus the above expression is the sum of the terms: first, second, third and eighth. Adding the remaining terms and taking common $[k-j]_q$ we get

$$\mathbf{D} = \eta s \gamma_{n,k-1,k-1-j} \, q^{n-1}\{[n-j]_q q^{j+1} + [n+1-k+j]_q$$
$$+ [k-1-j]_q q^{n+1-k+j} + [j]_q q - [n]_q - [n]_q q\}.$$

Thus (30) implies that the sum of the first, fourth and sixth term is zero as well as the sum of the second third and fifth terms. Consequently $\mathbf{D} = 0$.

**Acknowledgment.** The authors thank M. Ismail for a helpful discussion.

W. Bryc
Department of Mathematics
University of Cincinnati
PO Box 210025
Cincinnati, Ohio 45221-0025
USA
E-mail: Wlodzimierz.Bryc@UC.edu
URL: http://math.uc.edu/˜brycw

W. Matysiak
Faculty of Mathematics and Information Science
Warsaw University of Technology
Pl. Politechniki 1
00-661 Warszawa
Poland
and
Department of Mathematics
University of Cincinnati
PO Box 210025
Cincinnati, Ohio 45221-0025
USA
E-mail: matysiak@mini.pw.edu.pl
URL: http://www.mini.pw.edu.pl/˜matysiak

J. Wesołowski
Faculty of Mathematics and Information Science
Warsaw University of Technology
Pl. Politechniki 1
00-661 Warszawa
Poland
E-mail: wesolo@mini.pw.edu.pl
URL: http://www.mini.pw.edu.pl/tiki-index.php?page=prac_wesolowski_jacek